\documentclass[a4paper, 11pt]{amsart}
\usepackage{graphicx}
\usepackage{amssymb}
\usepackage[margin=3cm]{geometry}
\usepackage{color}
\usepackage{enumerate}
\usepackage[T1]{fontenc}
\usepackage[utf8]{inputenc}
\usepackage{mathtools}
\usepackage{comment}
\usepackage{float}
\usepackage{color}
\usepackage{mathrsfs}

\newtheorem{thm}{Theorem}[section]
\newtheorem{cor}[thm]{Corollary}

\theoremstyle{definition}
\newtheorem{defn}[thm]{Definition}
\theoremstyle{remark}

\numberwithin{equation}{section}

\begin{document}
\title[]{On Trace Theorems for Sobolev Spaces}

\author{Pier Domenico Lamberti}%
\address{Universit\`a degli Studi di Padova, Dipartimento di Matematica ``Tullio Levi-Civita'', Via Trieste 63, 35121 Padova, Italy}%
\email{lamberti@math.unipd.it}%
\author{Luigi Provenzano}%
\address{Universit\`a degli Studi di Padova, Dipartimento di Matematica ``Tullio Levi-Civita'', Via Trieste 63, 35121 Padova, Italy}%
\email{luigi.provenzano@math.unipd.it}%
\date{\today}%
\begin{abstract}
We survey a few trace theorems for Sobolev spaces on $N$-dimensional Euclidean domains. We include known results on linear subspaces, in particular  hyperspaces,  and smooth boundaries, as well as less known   results for Lipschitz boundaries, including  Besov's Theorem  and  other characterizations of traces on planar domains, polygons in particular, in the spirit of the work of P.~Grisvard. Finally, we present a recent approach, originally developed by G.~Auchmuty in the case of the Sobolev space $H^1(\Omega)$ on  a Lipschitz domain $\Omega$, and which we have further developed  for the trace spaces of $H^k(\Omega)$, $k\geq 2$, by using  Fourier expansions  associated with the eigenfunctions of new multi-parameter polyharmonic Steklov problems.
\end{abstract}

\maketitle




\section{Introduction}

The purpose of the present paper is twofold. First, we survey a few known  and less known results on the traces of functions of the Sobolev space $W^{k,p}(\Omega)$, $k\in\mathbb N$ and $1\leq p<\infty$, when $\Omega$ is a smooth or non-smooth  open set in $\mathbb R^N$,   $N\geq 2$. Second, for  bounded open sets $\Omega$ with Lipschitz boundaries, that is, open sets of class $C^{0,1}$,  we provide  an  explicit description  of the trace spaces of   $W^{k,2}(\Omega)$, which extends to arbitrary values of $k\geq 2$ the  results proved  in \cite{lamberti_provenzano_traces}  for $k=2$ and  based on new Steklov-type problems for polyharmonic operators.   In the sequel, the spaces $W^{k,2}(\Omega)$ will also be denoted by $H^k(\Omega)$.
 
Recall that if  $\Omega$ is a bounded  open set in $\mathbb R^N$ with Lipschitz boundary,  there exists a linear and continuous operator $\Gamma$ from $W^{k,p}(\Omega)$ to $(L^p(\partial\Omega))^k$ defined by $\Gamma(u)=(\gamma_0(u),...,\gamma_{k-1}(u))$, where $\gamma_0(u)$ is the {\it trace} of $u$ and $\gamma_j(u)$ is the $j$-th normal derivative of $u$ for $j=1,...,k-1$. In particular, for $u\in C^{k}(\overline\Omega)$, we have $\gamma_0(u)=u_{|_{\partial\Omega}}$ and $\gamma_j(u)=\frac{\partial^ju}{\partial\nu^j}$ for all $j=1,...,k-1$, where $\nu$ denotes the outer unit normal to $\partial\Omega$. The vector $\Gamma (u)$ is called the {\it total trace} of $u$.

Important problems in the theory of Sobolev spaces include the description of the {\it trace spaces} $\gamma_j(W^{k,p}(\Omega))$ for $j=0,...,k-1$, and the description of the {\it total trace space} $\Gamma(W^{k,p}(\Omega))$. 
 From a historical point of view, this problem finds its roots back at least  in 1906  with the publication  of the  paper \cite{hadamard} where 
J. Hadamard provided his famous counterexample which pointed out the need to clarify which conditions on the datum $g$ guarantee that the solution $v$ to the Dirichlet problem 
$$
\left\{ 
\begin{array}{ll}\Delta v=0,& {\rm in}\ \Omega,\\
v=g,& {\rm on }\ \partial \Omega,
\end{array}
\right.
$$ 
has square summable gradient.  
 Note that, in the framework of Sobolev spaces,  this  problem  can be reformulated as the problem of finding necessary and sufficient conditions on $g$ such that $g=\gamma_0(u)$ for some $u\in H^1(\Omega)$.  

Conclusive results are available for smooth domains, and are nowadays classical. The standard approach consists in flattening the boundary of $\Omega$ by means of suitable local diffeomorphisms. Hence the problem is recast to that of describing the trace spaces of $W^{k,p}(\mathbb R^N)$ on $N-1$-dimensional hyperplanes which can be identified with $\mathbb R^{N-1}$. A first classical method for describing the trace spaces of $W^{k,p}(\mathbb R^N)$ on $\mathbb R^{N-1}$ in the case $p=2$ is via Fourier Transform, see  e.g.,   \cite{lions_jl,necas}  and  Theorem~\ref{trace_thm_fourier}.  This method was already used in \cite{prodi}.   If $p\ne 2$ this approach is no more applicable. In this case, the description of the trace spaces relies on a method originally developed by E. Gagliardo in \cite{gagliardo} for the case $k=1$  in order to extend results obtained in \cite{aron,  nikolskii, slobo}  for $p=2$. See Theorems~\ref{trace_thm_gagliardo}, \ref{trace_thm_gagliardo0}. This method involves the use of Besov spaces $B^{s}_{p}(\mathbb R^{N-1})$, $s>0$, which are equivalent for non-integer $s$ to the fractional Sobolev spaces $W^{s,p}(\mathbb R^{N-1})$ appearing in \cite{gagliardo}. Applications of this method for $k\geq 2$ can be found in \cite{burenkov,grisvard,necas}. If the domain is sufficiently smooth, the definition of Besov spaces can be transplanted  from $\mathbb R^{N-1}$ to $\partial\Omega$, providing  well-defined function spaces $B^s_p(\partial\Omega)$ at the boundary of $\Omega$. As a matter of fact, it turns out that $\gamma_j(W^{k,p}(\Omega))=B_p^{k-j-1/p}(\partial\Omega)$ for all $j=0,...,k-1$ and $\Gamma(W^{k,p}(\Omega))=\prod_{j=0}^{k-1}B_p^{k-j-1/p}(\partial\Omega)$ for $p\ne 1$, see Theorem \ref{trace_thm_gagliardo}. 

However, when $\Omega$ is an arbitrary bounded open set with  Lipschitz boundary, there is no such  simple description and not many results are available in the literature. We collect some of them in this paper. First of all, we note that a complete description of the traces of all derivatives up to the order $k-1$ of a function $u\in W^{k,p}(\Omega)$ is due to O. Besov who provided an explicit but quite technical representation theorem, see ~\cite{besov19721, besov19722}, see also \cite{besov2}  and Theorem~\ref{besov_thm}. Simpler descriptions are not available with the exception of a few special cases. For example, when $\Omega$ is a curvilinear polygon in $\mathbb R^2$ with smooth edges, the trace spaces are described by using the classical trace spaces on each side of the polygon complemented with suitable compatibility conditions at the vertexes. This approach has been discussed by P. Grisvard in  the celebrated monographs \cite{grisvard,grisvard2}, see  Theorems \ref{polygon_1} and \ref{polygon_2}.    In \cite{grisvard2} one can also find a few related results on three-dimensional polyhedra. For more general planar domains and $p=2$, $k=2$, another description, given in terms of simple compatibility conditions is provided in \cite{geymonat_airy}, see  Theorem \ref{geymonat_1}. Theorem \ref{geymonat_1} is extended to the case $p\ne 2$ in \cite{duran_muschietti} and to the case $N=3$, $p\ne 2$ in \cite{geymonat}, see Theorem \ref{geymonat_3}. Moreover, necessary conditions for traces of functions in $W^{k,p}(\Omega)$ for all $k\geq 2$ are given in \cite{geymonat}.

We note that  our list  of results  is not exhaustive and we refer to the recent monograph  \cite{mazya_polyhedral} which  treats the trace problem in presence of corner or conical singularities in $\mathbb R^3$, as well as further results on $N$-dimensional polyhedra. We also quote  the fundamental paper ~\cite{kon} by  V.~Kondrat'ev for a pioneering work in this type of problems.  

A recent approach to trace spaces has been developed for $k=1$, $p=2$  by G. Auchmuty in \cite{auchmuty_steklov} 
where an alternative equivalent description of the trace space $\gamma_0(H^{1}(\Omega))$ is given  in terms of Fourier series associated with the eigenfunctions of the classical Steklov problem  \eqref{Steklov_classic} for the Laplace operator (see also \cite{touhami} for related results). This method has been employed  in \cite{lamberti_provenzano_traces}  for the case $k=2$, where new families of multi-parameter biharmonic Steklov  problems have been introduced with the specific purpose of describing the traces of functions in $H^2(\Omega)$. We emphasize the fact  that the  description of the trace spaces $\gamma_0(H^2(\Omega))$ and $\gamma_1(H^2(\Omega))$ and of the total trace space $\Gamma(H^2(\Omega))$ given 
 in \cite{lamberti_provenzano_traces}  is valid for arbitrary bounded open sets with Lipschitz boundaries in $\mathbb R^N$ and any  $N\geq 2$. 
 
In the present paper we generalize the result of \cite{lamberti_provenzano_traces} to the case $k\geq 2$. Following \cite{lamberti_provenzano_traces}, we provide decompositions of the space $H^k(\Omega)$ of the form $H^k(\Omega)=H^k_{0,\ell}(\Omega)+\mathcal H^k_{\ell}(\Omega)$ for all $\ell=0,...,k-1$. The spaces $H^k_{0,\ell}(\Omega)$ are the subspaces of $H^k(\Omega)$ of those functions $u$ such that $\gamma_{\ell}(u)=0$. The spaces $\mathcal H^k_{\ell}(\Omega)$ are associated with families of polyharmonic Steklov problems which we introduce in \eqref{multi-Steklov}, and admit Fourier bases of Steklov eigenfunctions, see Theorem \ref{steklov_thm_k}. Under the sole assumption that $\Omega$ is of class $C^{0,1}$ we use those bases to define in a natural way $k$ spaces at the boundary which we denote by $H^{k-\ell-1/2}_{\mathcal A}(\partial\Omega)$ for all $\ell=0,...,k-1$ (see \eqref{multi-trace-spaces} for precise definition) and prove that $\gamma_{\ell}(H^k(\Omega))=\gamma_{\ell}(\mathcal H^k_{\ell}(\Omega))=H^{k-\ell-1/2}_{\mathcal A}(\partial\Omega)$ for all $\ell=1,...,k-1$, see Theorem \ref{thm_single_traces}. It follows in particular that, if one wishes to define the space $H^{k-1/2}(\partial\Omega)$ as $\gamma_0(H^k(\Omega))$, our result gives an explicit description of  $H^{k-1/2}(\partial\Omega)$.

It turns out that the analysis of problems \eqref{multi-Steklov} provides further information on the total trace $\Gamma(H^k(\Omega))$. In particular $\Gamma(H^k(\Omega))\subset\prod_{\ell=0}^{k-1}H^{k-\ell-1/2}_{\mathcal A}(\partial\Omega)$. The inequality is in general strict if we assume that $\Omega$ is only of class $C^{0,1}$, see \cite{lamberti_provenzano_traces} for the case $k=2$. Moreover, we provide sufficient conditions for a $k$-tuple in  $\prod_{\ell=0}^{k-1}H^{k-\ell-1/2}_{\mathcal A}(\partial\Omega)$ to belong to $\Gamma(H^k(\Omega))$, see Theorem \ref{comp_thm}, see also Corollary \ref{comp_thm_2} for $k=2$. If $\Omega$ is smooth, we recover the classical result, namely $\Gamma(H^k(\Omega))=\prod_{\ell=0}^{k-1}H^{k-\ell-1/2}_{\mathcal A}(\partial\Omega)$ and in particular the spaces $H^{k-\ell-1/2}_{\mathcal A}(\partial\Omega)$ coincide with the classical trace spaces.

This paper is organized as follows. In Section \ref{pre} we introduce some notation and discuss a few preliminary results on the notion of trace. In Section \ref{classic} we review the classical trace theorems on smooth domains. In particular, in Subsection \ref{classic_1} we present the approach via  Fourier Transform, while in Subsection \ref{sub_gagliardo} we discuss Gagliardo's approach and the corresponding use of Besov spaces.  In Section \ref{nonsmooth} we review a few results on Lipschitz domains. In particular, in Subsection \ref{nonsmooth_1} we state Besov's Theorem. In Subsection \ref{nonsmooth_2} we collect a few results on curvilinear polygons in $\mathbb R^2$ and more general planar domains. In Subsection \ref{nonsmooth_3} we briefly describe the approach of  G. Auchmuty for the trace space of $H^1(\Omega)$ and we announce our results for  the general case of $H^k(\Omega)$ with $k\geq 2$ based on new Steklov problems for polyharmonic operators.

\section{Preliminaries on the notion of trace}\label{pre}

Let $\Omega$ be an open set in  $\mathbb R^N$, $1\leq p< \infty$ and $k\in\mathbb N$.  By $W^{k,p}(\Omega)$ we denote the Sobolev space  of functions $u\in L^p(\Omega)$ with all weak  derivatives of any order up to $k$ in $L^p(\Omega)$.  The space $W^{k,p}(\Omega)$ is  endowed with the norm
\begin{equation*}
\|u\|_{W^{k,p}(\Omega)}:=\biggl(\sum_{|\alpha|\leq k}\|D^{\alpha}u\|_{L^p(\Omega)}^p\biggr)^{\frac{1}{p}}.
\end{equation*}
We denote by $W^{k,p}_0(\Omega)$ the closure of $C^{\infty}_c(\Omega)$ with respect to $\|\cdot\|_{W^{k,p}(\Omega)}$, where $C^{\infty}_c(\Omega)$ is the space of functions in $C^{\infty}(\Omega)$ with compact support  in $\Omega$.

\subsection{Traces of functions on a $N-1$ dimensional subspace}\label{sub_traces}

Following \cite{burenkov}, we describe here a general explicit  definition of the  trace $\mathcal T(u)$ of a function $u\in L^1_{loc}(\mathbb R^N)$ on a $N-1$-dimensional subspace of $\mathbb R^N$, say $\mathbb R^{N-1}$, which will be later  applied to functions $u\in W^{k,p}(\mathbb  R^N) $. To do so, we represent each point $x$ of $\mathbb R^N$ as $x=(x',x_N)$, with $x'=(x_1,...,x_{N-1})\in\mathbb R^{N-1}$. Note that the space $\mathbb R^{N-1}$ is identified with the subspace of $\mathbb R^N$ of all points $(x',x_N)\in\mathbb R^N$ such that $x_N=0$.

For a continuous function $u$, the trace of $u$ on $\mathbb R^{N-1}$ is defined in a natural way as the restriction of $u$ on $\mathbb R^{N-1}$, namely $\mathcal T(u):=u_{|_{\mathbb R^{N-1}}}=u(x',0)$, $x'\in\mathbb R^{N-1}$. On the other hand,   this definition does not make sense for an arbitrary  function $u\in L^1_{loc}(\mathbb R^N)$  since it is  defined up to a set of zero Lebesgue measure.  As in \cite{burenkov},  we state the requirements which need to be fulfilled by a meaningful definition of the trace $g$ of a function $u\in L^1_{loc}(\mathbb R^N)$:
\begin{enumerate}[i)]
\item $g\in L^1_{loc}(\mathbb R^{N-1})$.
\item $f,g\in L^1_{loc}(\mathbb R^{N-1})$ are traces of a function $u\in L^1_{loc}(\mathbb R^N)$ if and only if they are equivalent on $\mathbb R^{N-1}$.
\item If $g\in L^1_{loc}(\mathbb R^{N-1})$ is the trace of $u\in L^1_{loc}(\mathbb R^N)$ and $v$ is equivalent to $u$ on $\mathbb R^N$, then $g$ is also the trace of $v$.
\item If $u$ is continuous, then $u(x',0)$ is the trace of $u$. 
\end{enumerate}
The following definition, given  in \cite[\S 5]{burenkov}, fulfills all the requirements above.
\begin{defn}
\label{trace}
Let $u\in L^1_{loc}(\mathbb R^N)$ and $g\in L^1_{loc}(\mathbb R^{N-1})$. We say that  the function $g$ is a trace of the function $u$ and we write $g=\mathcal T(u)$, if there exists a function $h\in L^1_{loc}(\mathbb R^N)$ equivalent to $u$ on $\mathbb R^N$ and
$$
h(\cdot,x_N)\rightarrow g(\cdot){\ \rm in\ }L^1_{loc}(\mathbb R^{N-1}){\ \rm as\ }x_N\rightarrow 0.
$$
\end{defn}
Other equivalent definitions of traces can be found e.g., in \cite{necas}. We have the following theorem on the existence of the traces of functions and of their derivatives.
\begin{thm}\label{trace_existence}
Let $k\in\mathbb N$ and $1\leq p< \infty$. Then the  traces $\mathcal T(D^{\alpha}u)$ on $\mathbb R^{N-1}$ of all weak partial derivatives $D^{\alpha}u$ with $|\alpha|\leq k-1$ exist and belong to $L^p(\mathbb R^{N-1})$ for all $u\in W^{k,p}(\mathbb R^N)$.
\end{thm}
Theorem~\ref{trace_existence} is usually proved  by establishing that the restriction of smooth functions to $\mathbb R^{N-1}$ defines a linear continuous operator which admits a unique  extension to the whole of $W^{k,p}(\mathbb R^N)$ which satisfies the requirements of  Definition~\ref{trace}.
 

Of particular interest is the description of the {\it total trace} $\Gamma(u)$ of a function $u\in W^{k,p}(\mathbb R^N)$ on $\mathbb R^{N-1}$, which is defined as the $k$-tuple
\begin{equation}\label{total_trace_flat}
\Gamma(u):=\left(\gamma_0(u),\gamma_1(u),...,\gamma_{k-1}(u)\right),
\end{equation}
where $\gamma_0(u)=\mathcal T(u)$ and $\gamma_j(u)=\mathcal T(\frac{\partial^j u}{\partial x_N^j})$ for $j=1,...,k-1$. We note that $\mathcal T(D^{\alpha}u)=D^{\alpha'}\mathcal T (\frac{\partial^{\alpha_{N}}u}{\partial x_N^{\alpha_{N}}})$, where $\alpha=(\alpha',\alpha_N)$ with $\alpha'\in\mathbb N_0^{N-1}$. This motivates the fact that we are interested only in the weak derivatives in the direction $x_N$. The following theorem holds. 

\begin{thm}
The map $\Gamma:W^{k,p}(\mathbb R^N)\rightarrow (L^p(\mathbb R^{N-1}))^k$ is a linear continuous  operator such that 
$$
\Gamma(u)=\biggl(u_{|_{\mathbb R^{N-1}}}, \frac{\partial u}{\partial {x_N}}_{|_{\mathbb R^{N-1}}},..., \frac{\partial^{k-1} u}{\partial {x_N^{k-1}}}_{|_{\mathbb R^{N-1}}}\biggr)
$$
for all $u\in C^{k}_c(\mathbb R^N)$.
\end{thm}

\subsection{Traces of functions on the boundary of an open set}\label{sub_traces_dom}

We recall now the notion of trace when $\mathbb R^N$ is replaced by $\Omega$ and $\mathbb R^{N-1}$ is replaced by $\partial\Omega$, where $\Omega\subset\mathbb R^N$ is a bounded domain, i.e., a bounded connected open set, and $\partial\Omega$ is its boundary. To do so, we need to describe suitable classes of domains. We recall the following definition where we use  the classical Schauder norm $\|\cdot \|_{C^{k,\gamma}}$ with the understanding that for $\gamma =0$ it coincides with the usual $\|\cdot \|_{C^{k}}$. 

\begin{defn}
\label{open_sets}
Let $\Omega\subset\mathbb R^N$ be a bounded domain. We say that $\Omega$ is a bounded domain of class $C^{k,\gamma}$ for some $k\in\mathbb N_0$ and $0\le \gamma\leq 1$ if there exist $\rho>0$, $s,s'\in\mathbb N$ with $s\leq s'$, a family $\left\{V_j\right\}_{j=1}^s$ of bounded open cuboids of the form $V_j=(a_{1j},b_{1j})\times\cdots\times (a_{Nj},b_{Nj})$ and a family $\left\{R_j\right\}_{j=1}^s$ of isometries in $\mathbb R^N$ such that
\begin{enumerate}[i)]
\item $\Omega\subset\cup_{j=1}^sV_j^{\rho}$ and $\Omega\cap V_j^{\rho}\ne\emptyset$ for all $j=1,...,s$, where  $V_j^{\rho}$ is defined by $V_j^{\rho}=\left\{x\in V_j:{\rm dist}(x,\partial V_j)>\rho\right\}$;
\item $\partial\Omega\cap V_j\ne\emptyset$ for $j=1,...,s'$ and $\partial\Omega\cap V_j=\emptyset$ for $j=s'+1,...,s$;
\item for $j=1,...,s$
$$
R_j(V_j)=\left\{x\in\mathbb R^N:a_{ij}<x_i<b_{ij}, i=1,...,N\right\}
$$
and
$$
R_j(\Omega\cap V_j)=\left\{x\in\mathbb R^N:a_{Nj}<x_N<\varphi_j(x'), x'\in W_j\right\},
$$
where $x'=(x_1,...,x_{N-1})$, $W_j=\left\{x'\in\mathbb R^{N-1}:a_{ij}<x_i<b_{ij}, i=1,...,N-1\right\}$ and $\varphi_j$ are functions of class $C^{k,\gamma}$ defined on $\overline W_j$ (it is meant that for $j=s'+1,...,s$ then $\varphi_j(x')=b_{Nj}$ for all $x'\in\overline W_j$) such that $\|D^{\alpha}\varphi_j\|_{C^{k,\gamma}(\overline W_j)}\leq M$ for all $|\alpha|\leq k$.  Moreover, for $j=1,...,s'$ it holds
$$
a_{Nj}+\rho\leq\varphi_j(x')\leq b_{Nj}-\rho.
$$
\end{enumerate}
We say that $\Omega$ is a bounded domain of class $C^{k}$ if it is of class $C^{k,\gamma}$  with $\gamma =0$.
\end{defn}

Assume now that $\Omega$ is a bounded domain of class $C^{0,1}$.
For a fixed $j\leq s'$ as in Definition \ref{open_sets}, consider the set $R_j(\partial\Omega\cap V_j)$. This set has the form $\left\{x\in\mathbb R^N:x_N=\varphi_j(x'), x'\in W_j\right\}$. Note that  by applying to the set $R_j(\Omega\cap V_j)$ the map $\Phi_j$   defined by $\Phi_j(x',x_N)=(x',x_N-\varphi_j(x'))$, we have that $\Phi_j\circ R_j (\partial\Omega\cap V_j)=\left\{x\in\mathbb R^N:x_N=0, x'\in W_j\right\}$. Thus we can give the following. 

\begin{defn}\label{partial_trace}
We say that $g$ is the trace of a function $u\in L^1(\Omega\cap V_j)$ on $\partial\Omega\cap V_j$ if $g\circ R_j^{(-1)}\circ\Phi_j^{(-1)}$ is the trace of $u\circ R_j^{(-1)}\circ\Phi_j^{(-1)}$ on $\mathbb R^{N-1}$ in the sense of Definition \ref{trace}.
\end{defn}

By using a suitable partition of unity, we can give a definition of trace of a function   $u\in L^1(\Omega)$ on $\partial\Omega$. Indeed, for all $j=1,...,s$ there exist a partition of unity given by functions $\psi_j\in C^{\infty}_c(\mathbb R^N)$ such that $|\psi_j(x)|\leq 1$ for all $x\in\mathbb R^N$, ${\rm supp}\,\psi_j\subset V_j$ for all $j=1,...,s$, $\sum_{j=1}^s\psi_j=1$ on $\Omega$. 
Thus, a function $u$ can be decomposed  as $\sum_{j=1}^su_j\psi_j$ and  its trace on $\partial \Omega$ can be defined by means of the following

\begin{defn}
 Assume that a function $u\in L^1(\Omega)$ is written in the form  $u=\sum_{j=1}^su_j$ where ${\rm supp}\,u_j\in V_j$ and $u_j\in L^1(\Omega\cap V_j)$. If the functions $g_j$ are traces of the functions $u_j$ on $\partial\Omega\cap V_j$ in the sense of Definition \ref{partial_trace}, then the function $g=\sum_{j=1}^sg_j$ is said to be the trace of the function $u$ on $\partial\Omega$, and we write $g=\mathcal T(u)$.
\end{defn}
We note that the previous definitions are well-posed and we refer to \cite[\S 2]{burenkov}  for more details. 

We  introduce now the total trace operator on $\Omega$ in analogy with Definition \ref{total_trace_flat}. We recall that if $\Omega$ is of class $C^{0,1}$ then a outer unit normal $\nu$ is defined almost everywhere on $\partial\Omega$. (Note that by using a suitable partition of unity as above, it is possible to define a $L^{\infty}$ vector field on $\overline\Omega$ which coincides almost everywhere with the normal vector field on $\partial\Omega$, see e.g., \cite[\S 1.5]{grisvard}.) 
 For a function $u\in W^{k,p}(\Omega)$ we define the total trace as
\begin{equation}\label{total_trace}
\Gamma(u):=\left(\gamma_0(u),\gamma_1(u),...,\gamma_{k-1}(u)\right),
\end{equation}
where $\gamma_0(u)=\mathcal T(u)$ and $\gamma_j(u)=\sum_{|\alpha|=j}\frac{j!}{\alpha!}\mathcal T(D^{\alpha}u)\nu^{\alpha}$ for $j=1,...,k-1$. For simplicity, we will often write, with abuse of notation, $\gamma_j(u)=\frac{\partial^j u}{\partial \nu^j}$ for $j=1,...,k-1$.

 We have the following.
\begin{thm}\label{total_trace_0}
Let $\Omega$ be a bounded domain in $\mathbb R^N$ of class $C^{0,1}$. Then the map $\Gamma:W^{k,p}(\Omega)\rightarrow (L^p(\partial\Omega))^k$ is a bounded linear operator such that 
$$
\Gamma(u)=\biggl(u_{|_{\partial\Omega}}, \frac{\partial u}{\partial {\nu}}_{|_{\partial\Omega}},..., \frac{\partial^{k-1} u}{\partial {\nu^{k-1}}}_{|_{\partial\Omega}}\biggr)\, ,
$$
for all $u\in C^{k}(\overline\Omega)$.
\end{thm}
We refer e.g., to \cite{necas} for a proof of Theorem \ref{total_trace_0}.

We conclude this subsection by recalling a characterization of the spaces $W^{k,p}_0(\Omega)$ by means of the corresponding traces.
\begin{thm}
Let $\Omega$ be a bounded domain  in $\mathbb R^N$ of class $C^{0,1}$, $1\leq p<\infty$ and $u\in W^{k,p}(\Omega)$. Then $\mathcal T(D^{\alpha}u)=0$ for all $|\alpha|\leq k-1$ if and only if $u\in W^{k,p}_0(\Omega)$.
\end{thm}
In particular, if $\Omega$ is of class $C^{0,1}$ then $u\in W^{1,p}_0(\Omega)$ if and only if $u\in W^{1,p}(\Omega)$ and $\mathcal T(u)=0$. Moreover, $u\in W^{2,p}_0(\Omega)$ if and only if $u\in W^{2,p}(\Omega)$ and $\Gamma(u)=0$ when $\Omega$ is of class $C^{0,1}$ (see \cite[Thm. 4.12]{necas}).  More generally, if the domain $\Omega$ is sufficiently regular, the space $W^{k,p}_0(\Omega)$ can be characterized by means of $\Gamma$. Namely, we have the following
\begin{thm}
Let $\Omega$ be a bounded domain  in $\mathbb R^N$ of class $C^{k,1}$, $1\leq p <\infty$ and $u\in W^{k,p}(\Omega)$. Then $\Gamma(u)=0$ if and only if $u\in W^{k,p}_0(\Omega)$.
\end{thm}
We refer to \cite[Thm. 1.5.1.5]{grisvard} or \cite[Thm. 4.13]{necas} for the proof.

\section{Classical Trace Theorems: smooth case}\label{classic}
This section is devoted to a short review of classical trace theorems. The focus is on the description of the {\it total trace space} $\Gamma(W^{k,p}(\Omega))$ defined by
$$
\Gamma(W^{k,p}(\Omega)):=\left\{\Gamma(u):u\in W^{k,p}(\Omega)\right\}.
$$

\subsection{Trace spaces via  Fourier Transform}\label{classic_1}

As is  customary, we  denote the Sobolev spaces $W^{k,2}(\Omega)$ and $W^{k,2}_0(\Omega)$ also by  $H^k(\Omega)$ and $H^k_0(\Omega)$, respectively.  

Recall that when $\Omega=\mathbb R^N$, the spaces $H^k(\Omega)$  can be equivalently defined  via Fourier Transform since $H^k(\mathbb R^N)$ is the space 
of functions $u\in  L^2( \mathbb R^N)$ such that 
\begin{equation}\label{sobolev_fourier}
\left(\int_{\mathbb R^N}(1+|\xi|^2)^k|\hat u(\xi)|^2d\xi\right)^{\frac{1}{2}}
\end{equation}
is finite. 
Here $\hat u $ denotes the Fourier Transform $\mathscr F[u]$  of a function $u\in L^2(\mathbb R^N)$ defined by 
\begin{equation*}
\hat u(\xi)=\mathscr F[u](\xi)=(2\pi)^{-\frac{N}{2}}\int_{\mathbb R^N}u(x)e^{-ix\cdot\xi}dx.
\end{equation*}

Recall that the left-hand side of \eqref{sobolev_fourier} defines a norm in $H^k(\mathbb R^N)$ equivalent to the standard one since
$$
\|u\|^2_{H^k(\mathbb R^N)}=\sum_{|\alpha|\leq k}\int_{\mathbb R^N}|D^{\alpha}u|^2dx=\int_{\mathbb R^N}\sum_{|\alpha|\leq k}|\xi^{\alpha}|^2|\hat u(\xi)|^2d\xi.
$$

The previous definitions extend to the case of  non-integer order of smoothness and allow to define the whole scale of spaces $H^s(\mathbb R^N)$, $s>0$ simply by replacing $k$ by $s$ in \eqref{sobolev_fourier}.

We have the following
\begin{thm}\label{trace_thm_fourier}
Let $k\in\mathbb N$. Then
$$
\Gamma( H^k(\mathbb R^N))=\prod_{j=0}^{k-1} H^{k-j-\frac{1}{2}}(\mathbb R^{N-1}).
$$
In particular, there exists $C>0$ such that 
$$
\|\Gamma(u)\|_{\prod_{j=0}^{k-1} H^{k-j-\frac{1}{2}}(\mathbb R^{N-1})}\leq C\|u\|_{H^k(\mathbb R^N)}\, ,
$$
for all $u \in H^k(\mathbb R^N)$. Moreover there exists a linear and continuous operator 
$$E:\prod_{j=0}^{k-1} H^{k-j-\frac{1}{2}}(\mathbb R^{N-1})\rightarrow  H^k(\mathbb R^N)$$
 such that if $u\in  H^k(\mathbb R^N)$, $u=Eg$ with $g\in \prod_{j=0}^{k-1} H^{k-j-\frac{1}{2}}(\mathbb R^{N-1})$, then $g=\Gamma(u)$.
\end{thm}
We note that  proving Trace Theorems consists of proving two statements: an embedding and an extension theorem.  In this case, the proof of the  embedding $\Gamma( H^k(\mathbb R^N))\subset\prod_{j=0}^{k-1} H^{k-j-1/2}(\mathbb R^{N-1})$ is straightforward. For example, for $\gamma_0(u)=\mathcal T (u)$ it is sufficient to write for a function $u\in C^{\infty}_c(\mathbb R^N)$, $\mathcal T u(x')=u(x',0)=\mathscr F^{(-1)}[\hat u](x',0)$. Fubini-Tonelli's Theorem, H\"older's inequality and standard manipulations allow to prove quite easily that $\|\mathcal T(u)\|_{ H^{k-\frac{1}{2}}(\mathbb R^{N-1})}\leq C \|u\|_{ H^k(\mathbb R^N)}$. The result is extended to $H^k(\mathbb R^N)$ by standard approximation.

As for the extension theorem, starting from an element $g\in\prod_{j=0}^{k-1} H^{k-j-1/2}(\mathbb R^{N-1})$, one constructs  explicitly a function $u\in  H^k(\mathbb R^N)$ which turns out to have total trace $g$ on $\mathbb R^{N-1}$. Namely, if $g=(g_0,...,g_{k-1})\in\prod_{j=0}^{k-1} H^{k-j-1/2}(\mathbb R^{N-1})$ , we define
$$
u(x',x_N)=\mathscr F^{(-1)}_{\xi'}\left[\sum_{j=0}^{k-1}\frac{x_N^j}{j!}\mathscr F_{\xi'}[g_j](\xi')h(x_N\sqrt{1+|\xi'|^2})\right],
$$
for some $h\in C^{\infty}_c(\mathbb R)$, $0\leq h(t)\leq 1$ and $h(t)=1$ for $|t|\leq 1$. Here $\xi'$  is defined by  $\xi'=(\xi_1,...,\xi_{N-1})$ and the Fourier Transform $\mathscr F_{\xi'}$ and its inverse are taken with respect to the variable $\xi'\in\mathbb R^{N-1}$. We refer to \cite[\S 2.5]{necas} for the details of the proof of Theorem \ref{trace_thm_fourier}.

Note that the space $H^s(\mathbb R^N)$ is defined for any $s\in\mathbb R$ by replacing $k$ by $s$ in \eqref{sobolev_fourier}. It is then possible to extend the validity of  Theorem \ref{trace_thm_fourier} to the case  non-integer $s\in\mathbb R$, as long as $s>\frac{1}{2}$. In this case one sees that
$$
\Gamma( H^s(\mathbb R^N))=\prod_{j=0}^{[s-\frac{1}{2}]} H^{s-j-\frac{1}{2}}(\Omega),
$$
and the other statements of Theorem \ref{trace_thm_fourier} remain valid.

\subsection{Gagliardo's method and Besov spaces}\label{sub_gagliardo}
In the case of a bounded domain  $\Omega$ in $\mathbb R^N$ the trace spaces of $W^{k,p}(\Omega)$ can described by means of  Gagliardo-Slobodeckij norms  which can also be encoded in the  Besov  spaces $B_p^s(\partial\Omega)$.
We begin  by recalling the definition of the Besov spaces in $\mathbb R^N$.  For $\ell\in\mathbb N$ and $h\in\mathbb R^N$ we define the difference of order $\ell$ of a function $f$ with step $h$ as
$$
\Delta_h^{\ell}u(x):=\sum_{j=0}^{\ell}(-1)^{\ell-j}\binom{\ell}{j}u(x+jh).
$$
\begin{defn}
\label{besov}
Let $s>0$ and $1\leq p <\infty$. Let  $\sigma\in\mathbb N$, $\sigma>s$. A function $u\in L^1_{loc}(\mathbb R^N)$ belongs to the Besov  space $B_{p}^s(\mathbb R^N)$ if 
$$
\|u\|_{B_{p}^s(\mathbb R^N)}:=\|u\|_{L^p(\mathbb R^N)}+|u|_{B_{p}^s(\mathbb R^N)} <\infty\, ,
$$ 
where
$$
|u|_{B_{p}^s(\mathbb R^N)}:=\left(\int_{\mathbb R^N}{\frac{\|\Delta^{\sigma}_h u\|_{L^p(\mathbb R^N)}^p}{|h|^{sp+N}}dh}\right)^{\frac{1}{p}}\, .
$$
\end{defn}
We remark that Definition \ref{besov} does not depend on the choice of $\sigma\in\mathbb N$, $\sigma>s$,  see e.g., \cite[\S 5.3]{burenkov}.

We recall that when $s>0$ is not an integer number, the space $B_p^s(\mathbb R^N)$ coincides with the Gagliardo-Slobodeckij space $W^{s,p}(\mathbb R^N)$ which is defined as the space of functions in $W^{[s],p}(\mathbb R^N)$  such that 
$$
|D^{\alpha}u|_{W^{s-[s],p}(\mathbb R^N)}:=\int_{\mathbb R^N}\int_{\mathbb R^N}\frac{|D^{\alpha}u(x)-D^{\alpha}u(y)|^p}{|x-y|^{ p (s-[s])+N}}dxdy<\infty,
$$
for all $\alpha\in\mathbb N_0^N$ with  $|\alpha|=[s]$, where $[s]$ is  the integer part of $s$. The space $W^{s,p}(\mathbb R^N)$ is endowed with the norm   
$$
\|u\|_{W^{s,p}(\mathbb R^N)}:=\biggl(\|u\|_{W^{[s],p}(\mathbb R^N)}^p+\sum_{|\alpha|=[s]}\int_{\mathbb R^N}\int_{\mathbb R^N}\frac{|D^{\alpha}u(x)-D^{\alpha}u(y)|^p}{|x-y|^{p(s-[s])+N}}dxdy\biggr)^{\frac{1}{p}},
$$
and the quantity $|\cdot |_{W^{s-[s],p}(\mathbb R^N)}$ is often called  {\it Gagliardo semi-norm}. 
We refer to \cite{palatucci} for more information on fractional Sobolev spaces defined on more general open sets of $\mathbb R^N$. We also remark that for all $s>0$ and for $p=2$, the space $B_2^s(\mathbb R^N)$ coincides with the space $ H^s(\mathbb R^N)$ defined via Fourier Transform as in \eqref{sobolev_fourier} with $k$ replaced by $s$, and that the two corresponding norms are equivalent.

We now define the  Besov spaces  $B_p^s(\partial\Omega)$ on $\partial\Omega$. To do so, we use similar arguments and notation  as in Subsection \ref{sub_traces_dom}. 

\begin{defn}
\label{besov_2}
Let $k\in\mathbb N$, $1\leq p <\infty$ and let $\Omega$ be a bounded domain of class $C^{k}$. Let $s< k$. We say that $g\in B_p^s(\partial\Omega)$ if 
$$
\|g\|_{B_p^s(\partial\Omega)}:=\left(\sum_{j=1}^{s'}\|(g\psi_j)\circ R_j^{(-1)}\circ\Phi^{(-1)}\|_{B_p^s(\Phi_j\circ R_j (\partial\Omega\cap V_j))}\right)^{\frac{1}{p}}<\infty.
$$
\end{defn}
Definition \ref{besov_2} does not depend on the particular choice of the cuboids $V_j$ and of the partition of unity $\psi_j$, see \cite[Ch. 5, Rem. 19]{burenkov}. For more details on Definition \ref{besov_2} and for more information on Besov spaces on smooth boundaries, we refer to \cite[\S 5]{burenkov}. 
 We remark that the norm of $B^s_p(\partial\Omega)$ when $0<s<1$ and  $\Omega\subset\mathbb R^N$ is a bounded domain of class $C^{0,1}$ can be given either by using Definition \ref{besov_2} or equivalently by setting
\begin{equation}\label{gagliardo_norm}
\|u\|_{W^{s,p}(\partial\Omega)}:=\left(\|u\|_{L^p(\partial\Omega)}^p+\int_{\partial\Omega}\int_{\partial\Omega}\frac{|u(x)-u(y)|^p}{|x-y|^{sp+N-1}}d\sigma(x)d\sigma(y)\right)^{\frac{1}{p}},
\end{equation}
see e.g., to \cite{grisvard}. In fact, the norm \eqref{gagliardo_norm} is the one which originally appears in the paper of E. Gagliardo \cite{gagliardo} 
where the following theorem was proved for $k=1$.  For the proof of Theorem \ref{trace_thm_gagliardo} in the case $k\geq 2$, we refer to \cite{burenkov,necas} for $1<p<\infty$ and to \cite{burenkov} for $p=1$.

\begin{thm}\label{trace_thm_gagliardo}
Let $k\in\mathbb N$ and let $\Omega$ be a bounded domain in $\mathbb R^N$ of class $C^{k,1}$. Then
$$
\Gamma(W^{k,p}(\Omega))=\prod_{j=0}^{k-1}B^{k-j-\frac{1}{p}}_p(\partial\Omega),\ \ \ 1<p <\infty,
$$
and
$$
\Gamma(W^{k,1}(\Omega))=\prod_{j=0}^{k-2}B^{k-j-1}_1(\partial\Omega)\times L^1(\partial\Omega).
$$
Moreover, $\Gamma$ is a continuous operator between $W^{k,p}(\Omega)$ and the corresponding total trace space. 
\end{thm}

We note that for $k=1$ the regularity assumptions on $\Omega$ in the previous theorem can be relaxed. In fact, the original result proved by Gagliardo in \cite{gagliardo} for $k=1$ requires that $\Omega$ is of class $C^{0,1}$ and reads as follows

\begin{thm}\label{trace_thm_gagliardo0}
Let $\Omega$ be a bounded domain in $\mathbb R^N$ of class $C^{0,1}$. Then
$$
\gamma_0(W^{1,p}(\Omega))=B^{1-\frac{1}{p}}_p(\partial\Omega),\ \ \ 1<p <\infty,
$$
and
$$
\gamma_0(W^{1,1}(\Omega))= L^1(\partial\Omega).
$$
Moreover, $\gamma_0$ is a continuous operator between $W^{1,p}(\Omega)$ and the corresponding trace space. 
\end{thm}

It is interesting to observe that for $p=1$ the extension operator from $L^1(\partial\Omega)$ to $W^{1,1}(\Omega)$ provided in \cite{gagliardo} is nonlinear, see \cite{burenkov_goldman} for further results in this direction.

We observe that the previous theorem does not make any essential use of the Besov norm itself  since the Gagliardo-Slobodeckij norm is enough for stating  it. Indeed, since the codimension of the manifold $\partial \Omega$ is one, only fractional orders of smoothness are involved in the statement. However, Besov spaces play a crucial role in describing the trace spaces on sub-manifolds of codimension larger than one in which case 
integer orders of smoothness may appear. The following theorem is a special case of a result proved by O. Besov in \cite{besov0, besov1}
which provides  the original and main justification for the introduction 
of Besov   spaces in the literature.  To give an idea of this, we state  Besov's trace Theorem in its simplest form for the trace ${\mathcal {T}}_m(u)$ on the subspace ${\mathbb{R}}^m$ of $\mathbb R^N$ for a function $u$ defined in $\mathbb R^N$. Here the definition of ${\mathcal {T}}_m(u)$ can be given as in Definition~\ref{trace} with $m$ replacing $N-1$. Note that, apart from the special case when  $p=1$ and $k=N-m$,  a necessary and sufficient condition for the existence of   ${\mathcal {T}}_m(u)$ for $u\in W^{k,p}(\mathbb R^N)$ is that $pk>N-m$.

\begin{thm}\label{trace_m_besov}
Let $k, m\in\mathbb N$, $1\le m<N$ and  $1\leq p <\infty$.  Then
$$
{\mathcal T}_m(W^{k,p}(\mathbb R^N))=B^{k-\frac{N-m}{p}}_p(\mathbb R^m),\ \ \  {\rm if}\  pk>N-m,
$$
and
$$
{\mathcal T}_m(W^{N-m,1}(\mathbb R^N))=L^1(\mathbb R^m).
$$
Moreover, ${\mathcal T}_m$ is a continuous operator between $W^{k,p}(\mathbb R^N)$ and the corresponding  trace spaces. 
\end{thm}

We refer to \cite[Ch.~5]{burenkov} for a detailed proof.    Note that using Besov spaces is essential in the previous theorem  when $ k-\frac{N-m}{p}\in {\mathbb{N}}$.

\section{Trace Theorems: Lipschitz case}\label{nonsmooth}
When $\Omega$ is an arbitrary domain in $\mathbb R^N$ of class $C^{0,1}$ there is not a description of $\Gamma(W^{k,p}(\Omega))$ as simple  as the one  given by Theorem \ref{trace_thm_gagliardo}. Actually, in this case the definition of the spaces $B^s_p(\Omega)$ is problematic when $s>1$ and  not many results are available in the literature. We shall present a few of them in the present section.

\subsection{Besov's Theorem}\label{nonsmooth_1}
Let $\Omega$ be a bounded domain in $\mathbb R^N$ of class $C^{0,1}$.  Thus, there exist $s,s'\in\mathbb N$ , open cuboids $V_j$, isometries $R_j$ and Lipschitz functions $\varphi_j:W_j\rightarrow\mathbb R$ as in Definition \ref{open_sets}. Let us denote by $M_j$ the Lipschitz constant of $\varphi_j$, for all $j=1,...,s$. We introduce a few more definitions.
For $h>0$, we denote by $A_j^h$ the  cone
$$
A_j^h:=\left\{x=(x',x_N)\in\mathbb R^N: x_N>M_j|x'|,\  |x|<h\right\}
$$
for all $j=1,...,s$. We may assume, possibly choosing a different isometry $R_j$, that
$$
(\partial\Omega\cap V_j)+R_j^{(-1)}(A_j^h)\subset\Omega
$$
and that $((\partial\Omega\cap V_j)+R_j^{(-1)}(A_j^h))\cap\Omega$ coincides with a sufficiently small neighborhood of some point of $\partial\Omega$. We also set
$$
A_j:=\left\{x=(x',x_N)\in\mathbb R^N: |x_N|>(M_j+\varepsilon)|x'|,{\rm\ for\ some\ }\varepsilon>0\right\},
$$
$$
\partial V_j(x):=R_j^{(-1)}(R_j(\partial\Omega\cap V_j)\cap (x-A_j)),
$$
and
$$
\Omega_j^h:=\left\{x\in\Omega:{\rm dist}(x,\partial V_j\cap\Omega)<h\right\},
$$
for all $j=1,...,s$. We state the following theorem, which is proved in \cite[Ch.V,\,\S 20]{besov2}
\begin{thm}\label{besov_thm}
Let $\Omega$ be a bounded domain in $\mathbb R^N$ of class $C^{0,1}$ and let $k\in\mathbb N$. Then, for any $u\in W^{k,p}(\Omega)$ and any $\alpha\in\mathbb N_0^N$ with $|\alpha|\leq k-1$ there exist traces of the derivatives $D^{\alpha}u$ for which we have
\small
\begin{multline}\label{besov_condition}
\sum_{j=1}^{s'}\left(\sum_{|\alpha|\leq k-1}\int_{\partial\Omega\cap V_j}|D^{\alpha}u|^pd\sigma\right)^{\frac{1}{p}}\\
+\sum_{j=1}^{s'}\Bigg(\sum_{|\alpha|\leq k-1}\int_{\Omega_j^h}\int_{\partial V_j(x)}\int_{\partial V_j(x)}\left|\frac{D^{\alpha}_x\left(\sum_{|\beta|\leq k-1}(D^{\beta}_y u(y)(x-y)^{\beta}-D^{\beta}_zu(z)(x-z)^{\beta})/\beta!\right)}{{\rm dist}(x,\partial\Omega\cap V_j)^{k+2\frac{N-1}{p}-|\alpha|}}\right|^p\\
\times d\sigma(y)d\sigma(z)dx \Bigg)^{\frac{1}{p}}
\leq C\|u\|_{W^{k,p}(\Omega)},
\end{multline}
\normalsize
where $h>0$ is a sufficiently small number and the constant $C>0$ does not depend on $u$.

Conversely, suppose that a set $\left\{g_{\alpha}\right\}_{|\alpha|\leq k-1}$, $g_{\alpha}\in L^p(\partial\Omega)$ is such that the left-hand side of \eqref{besov_condition} with $D^{\beta}u$ replaced by $g_{\beta}$ is finite. Then, there exists $u\in W^{k,p}(\Omega)$ for which $\mathcal T (D^{\alpha}u)$ exist for all $|\alpha|\leq k-1$, $\mathcal T(D^{\alpha}u)=g_{\alpha}$ and $\|u\|_{W^{k,p}(\Omega)}$ is estimated by a constant independent on $u$ times the left-hand side of \eqref{besov_condition}, with $D^{\beta}u$ replaced by $g_{\beta}$.
\end{thm}

\subsection{Polygons and planar sets}\label{nonsmooth_2}
Simpler descriptions of the trace spaces of $W^{k,p}(\Omega)$  are available when $\Omega$ is a polygon in $\mathbb R^2$. We say that a bounded domain  $\Omega$ in $\mathbb R^2$ is a {\it curvilinear polygon} of class $C^{k,\gamma}$ for some $k\in\mathbb N_0$, $0\leq\gamma\leq 1$, if $\partial\Omega=\bigcup_{j=1}^n\overline{\Gamma_j}$, $\Gamma_i\cap\Gamma_{j}=\emptyset$ for $i\ne j$,  $\overline{\Gamma_j}\cap\overline{\Gamma_{j+1}}=\mathcal V_{j}$ for $j=1,...,n-1$, $\overline{\Gamma_n}\cap\overline{\Gamma_{1}}=\mathcal V_n$ and $\overline{\Gamma_i}\cap\overline{\Gamma_{j}}=\emptyset$ in the other cases, where $\Gamma_j\subset\mathbb R^2$ are curves of class $C^{k,\gamma}$ called {\it sides} of the polygon, and $\mathcal V_{j}\in\mathbb R^2$ are the {\it vertexes} of the polygon.

Theorem \ref{trace_thm_gagliardo} is easily seen to hold with $\Gamma$ replaced by $\Gamma_{|_{\Gamma_j}}$, i.e., the restriction to $\Gamma_j$ of the total trace operator $\Gamma$ defined by \eqref{total_trace}, and $\partial\Omega$ replaced by $\Gamma_j$, see e.g., \cite[\S 1.5]{grisvard}.

However, in many applications the knowledge of the traces on all the sides $\Gamma_j$ is not sufficient, and one looks for the image of $W^{k,p}(\Omega)$ on the whole of $\partial\Omega$ through the operator $\Gamma$. To do so, compatibility conditions at the vertexes are possibly needed, as highlighted e.g., in \cite{grisvard,grisvard2}. For $W^{1,p}(\Omega)$, $1<p<\infty$, we have the following theorem from \cite[Thm. 1.5.2.3]{grisvard}.
\begin{thm}\label{polygon_1}
Let $\Omega$ be a curvilinear polygon in $\mathbb R^2$ of class $C^{1}$ with boundary $\partial\Omega=\bigcup_{j=1}^n\overline{\Gamma_j}$ and let $1< p<\infty$. Then $\mathcal T$ is a linear and continuous mapping with continuous inverse from $W^{1,p}(\Omega)$ to the subspace of $\prod_{j=1}^nB^{1-1/p}_p(\Gamma_j)$ of functions  $(g_1,...,g_n)\in \prod_{j=1}^nB^{1-1/p}_p(\Gamma_j)$ satisfying:
\begin{enumerate}[i)]
\item no extra conditions, when $1<p<2$;
\item $g_j(\mathcal V_j)=g_{j+1}(\mathcal V_j)$ for all $j=1,...,n$, with the convention that $g_{n+1}=g_1$, when $2<p<\infty$;
\item $\int_0^{\delta_j}\frac{|g_{j+1}(x_j(\sigma))-g_{j}(x_j(-\sigma))|^2}{\sigma}d\sigma<\infty$ for all $j=1,...,n$, when $p=2$, where $x_j(\sigma)$ denotes the point on $\partial\Omega$ at arc-length distance $\sigma$ from $\mathcal V_j$, and $\delta_j>0$ is such that, when $|\sigma|\leq\delta_j$, then $x_j(\sigma)\in\Gamma_j$ if $\sigma>0$ and $x_j(\sigma)\in\Gamma_{j+1}$ if $\sigma<0$.
\end{enumerate}

\end{thm}

The result for $W^{k,p}(\Omega)$ for $1<p<\infty$ is  stated in \cite[Thm. 1.5.2.8]{grisvard} and reads as follows.

\begin{thm}\label{polygon_2}
Let $\Omega$ be a curvilinear polygon in $\mathbb R^2$ of class $C^{\infty}$ with boundary $\partial\Omega=\bigcup_{j=1}^n\overline{\Gamma_j}$ and let $1< p<\infty$. Then $\Gamma$ is a linear and continuous mapping with continuous inverse from $W^{k,p}(\Omega)$ to the subspace of
$$
\prod_{j=1}^n\prod_{i=0}^{k-1}B^{k-i-\frac{1}{p}}_p(\Gamma_j)
$$ 
given by  those elements  $((g^{(1)}_0,...,g^{(1)}_{k-1}),...,(g^{(j)}_0,...,g^{(j)}_{k-1}),...,(g^{(n)}_0,...,g^{(n)}_{k-1}))$ of  the space $\prod_{j=1}^n\prod_{i=0}^{k-1}B^{k-i-1/p}_p(\Gamma_j)$ defined by the following conditions: let $L$ be any linear differential operator with coefficients of class $C^{\infty}$ of order $m\leq k-\frac{2}{p}$; denote by $P_{j,i}$ the differential operator tangential to $\Gamma_j$ such that $L=\sum_{i\geq 0}P_{j,i}\frac{\partial^j}{\partial\nu_j^l}$, where $\nu_j$ denotes the outer unit normal to $\Gamma_j$; then
\begin{enumerate}[i)]
\item $\sum_{i\geq 0}(P_{j,i}g^{(j)}_i)(\mathcal V_j)=\sum_{i\geq 0}(P_{j,i}g^{(j+1)}_i)(\mathcal V_j)$ for $m<k-\frac{2}{p}$;
\item $\int_0^{\delta_j}\left|\sum_{i\geq 0}\left((P_{j,i}g_i^{(j)})(x_j(-\sigma))-(P_{j+1,i}g^{(j+1)}_i)(x_j(\sigma))\right)\right|^2\frac{d\sigma}{\sigma}<\infty$ for $p=2$ and  $m=k-1$.
\end{enumerate}

\end{thm}

For the proofs of Theorems \ref{polygon_1} and \ref{polygon_2} we refer to \cite{grisvard,grisvard2}. A few information on the compatibility conditions on the edges and the vertexes of three-dimensional polyhedra are available in \cite{grisvard2}. For a more detailed analysis on trace spaces on domains with corner and conical singularities in $\mathbb R^3$ and for trace spaces on $N$-dimensional polyhedra we also refer to the monograph \cite{mazya_polyhedral}.

A characterization of the range of $\Gamma(H^2(\Omega))=(\gamma_0(H^2(\Omega)),\gamma_1(H^2(\Omega)))$ in terms of compatibility conditions in the case when  $\Omega\subset\mathbb R^2$ is just of class $C^{0,1}$ has been given in \cite{geymonat_airy}. It is stated as follows.

\begin{thm}\label{geymonat_1}
Let $\Omega$ be a bounded domain in $\mathbb R^2$ of class $C^{0,1}$ and let $g_0\in H^1(\partial\Omega)$, $g_1\in L^2(\partial\Omega)$. Then there exists $u\in H^2(\Omega)$ such that $(g_0,g_1)=\Gamma(u)$ if and only if
\begin{equation}\label{comp_0}
(\partial_tg_0)\nu-g_1 t\in \big(B^{\frac{1}{2}}_2(\partial\Omega)\big)^2,
\end{equation}
where $t$ denotes the positively oriented unit tangent vector to $\partial\Omega$.
\end{thm}
Note that in the case of a smooth set, the vectors $\nu$ and $t$ are linearly independent at every point of $\partial\Omega$, thus one recovers the characterization given in Theorem \ref{trace_thm_gagliardo} with $k=2$, $p=2$, $N=2$. As pointed out in \cite{geymonat_airy}, the compatibility conditions \eqref{comp_0} are equivalent to those of Theorem \ref{polygon_2} when $k=2$, $p=2$. In \cite{duran_muschietti}, Theorem \ref{geymonat_1} is extended to the case $1<p<\infty$. An equivalent characterization of the range of $\Gamma$ on $W^{2,p}(\Omega)$ when $\Omega\subset\mathbb R^2$ is a Lipschitz domain is given in \cite{geymonat}. Namely, we have the following.
\begin{thm}\label{geymonat_2}
Let $\Omega $ be a bounded domain in $\mathbb R^2$ of class $C^{0,1}$ and let $g_0\in W^{1,p}(\partial\Omega)$, $g_1\in L^p(\partial\Omega)$. Then there exists $u\in W^{2,p}(\Omega)$ such that $(g_0,g_1)=\Gamma(u)$ if and only if
\begin{equation}\label{comp_01}
(\partial_tg_0)t+g_1\nu\in \big(B^{1-\frac{1}{p}}_p(\partial\Omega)\big)^2.
\end{equation}
\end{thm}
We refer to \cite{geymonat} for further discussions on Theorem \ref{geymonat_2}. Exploiting condition \eqref{comp_01} allows to provide a characterization of $\Gamma(W^{1,p}(\Omega))$ also for $N=3$. Indeed, the following theorem is proved in \cite{geymonat3}.

\begin{thm}\label{geymonat_3}
Let $\Omega $ be a bounded domain in $\mathbb R^3$ of class $C^{0,1}$ and let $g_0\in W^{1,p}(\partial\Omega)$, $g_1\in L^p(\partial\Omega)$. Then there exists $u\in W^{2,p}(\Omega)$ such that $(g_0,g_1)=\Gamma(u)$ if and only if
\begin{equation}\label{comp_02}
\nabla_{\partial\Omega}g_0+g_1\nu\in \big(B^{1-\frac{1}{p}}_p(\partial\Omega)\big)^3,
\end{equation}
where $\nabla_{\partial\Omega}g$ denotes the tangential gradient of $g$ on $\partial\Omega$.
\end{thm}

Necessary conditions for the traces of functions in $W^{k,p}(\Omega)$ for all $k\geq 2$, $1<p<\infty$ are given in \cite{geymonat}. These conditions are also sufficient for $p=2$ and $N=2$, thus recovering Theorem \ref{geymonat_1}. In \cite{geymonat} the authors present a general scheme to write necessary conditions which turn out to be of the form \eqref{comp_0}, \eqref{comp_01}, \eqref{comp_02} and write the condition for $k=3$ only. We refer to \cite[Thm. 3.4]{geymonat} for the precise statement.

\subsection{Auchmuty's method}\label{nonsmooth_3}
In this subsection we present a recent approach for describing the traces of functions in $H^k(\Omega)$ on the boundary of a Lipschitz domain $\Omega$ of $\mathbb R^N$. The trace spaces are defined by means of Fourier series associated with the eigenfunctions of families of Steklov-type problems for the polyharmonic operator $(-\Delta)^k$. The definitions of the trace spaces coincide with the classical ones when the domain is sufficiently smooth.

This approach has been developed by G. Auchmuty \cite{auchmuty_steklov_2,auchmuty_steklov} for the trace space $\Gamma(H^1(\Omega))=\mathcal T(H^1(\Omega))$, which is known to coincide with $H^{1/2}(\partial\Omega):=B^{1/2}_2(\partial\Omega)$, see Theorem \ref{trace_thm_gagliardo}. It has been recently extended in \cite{lamberti_provenzano_traces} in order to characterize $\Gamma(H^2(\Omega))$ when $\Omega\subset\mathbb R^N$ is a bounded Lipschitz domain. We will describe here how the results of \cite{lamberti_provenzano_traces} apply in general for any  $k\geq 2$.

\subsubsection{Case $k=1$} We find it convenient to describe the original method of Auchmuty for $k=1$  first. On a bounded domain $\Omega$ in $\mathbb R^N$ of class $C^{0,1}$ we consider the following variational eigenvalue problem
\begin{equation}\label{Steklov}
\int_{\Omega}\nabla u\cdot\nabla\varphi dx=\sigma\int_{\partial\Omega}u\varphi dx\,,\ \ \ \forall\varphi\in H^1(\Omega),
\end{equation}
in the unknowns $u\in H^1(\Omega)$ (the eigenfunction) and $\sigma\in\mathbb R$ (the eigenvalue). Problem \eqref{Steklov} is the weak formulation of the well-known Steklov eigenvalue problem, namely
\begin{equation}\label{Steklov_classic}
\begin{cases}
\Delta u=0, & {\rm in\ }\Omega,\\
\frac{\partial u}{\partial\nu}=\sigma u, & {\rm on\ }\partial\Omega.
\end{cases}
\end{equation}
We recall that a function $u\in H^1(\Omega)$ is called harmonic if
$$
\int_{\Omega}\nabla u\cdot\nabla\varphi dx=0\, ,
$$
for all $\varphi\in H^1_0(\Omega)$. We denote by $\mathcal H^1(\Omega)$  the space of all harmonic functions in $H^1(\Omega)$. We consider on $H^1(\Omega)$ the scalar product
\begin{equation}\label{equivalent_H1}
\langle u,v\rangle_{H^1_{\partial}(\Omega)}:=\int_{\Omega}\nabla u\cdot\nabla v dx+\int_{\partial\Omega}u vd\sigma\,,\ \ \ \forall u,v\in H^1(\Omega),
\end{equation}
which induces on $H^1(\Omega)$ the equivalent norm
\begin{equation}\label{equivalent_H1_norm}
\|u\|_{H^1_{\partial}(\Omega)}^2:=\int_{\Omega}|\nabla u|^2dx+\int_{\partial\Omega}u^2d\sigma\,,\ \ \ \forall u\in H^1(\Omega).
\end{equation}
Thus, we have the following decomposition of the space $H^1(\Omega)$
$$
H^1(\Omega)=H^1_0(\Omega)\oplus\mathcal H^1(\Omega),
$$
where the sum is orthogonal with respect to \eqref{equivalent_H1}.

We have the following theorem on the spectrum of problem \eqref{Steklov}, the proof of which can be found in \cite{auchmuty_steklov_2}.
\begin{thm}\label{steklov_thm_0}
Let $\Omega $ be a bounded domain in $\mathbb R^N$ of class $C^{0,1}$. The eigenvalues of problem \eqref{Steklov} have finite multiplicity and are given by a non-decreasing sequence of non-negative real numbers $\sigma_j$ defined by
\begin{equation*}
\sigma_j=\min_{\substack{U\subset H^1(\Omega)\setminus  H^1_0(\Omega) \\{\rm dim}\,U=j}}\max_{\substack{u\in U\\u\ne 0}}\frac{\int_{\Omega}|\nabla u|^2dx}{\int_{\partial\Omega}u^2d\sigma},
\end{equation*}
where each eigenvalue is repeated according to its multiplicity. The first eigenvalue $\sigma_1=0$ has multiplicity one and the corresponding eigenfunctions are the constant functions on $\Omega$. Moreover, there exists a Hilbert basis $\left\{u_j\right\}_{j=1}^{\infty}$ of $\mathcal H^1(\Omega)$ of eigenfunctions $u_j$.
 Finally, by normalizing the eigenfunctions $u_j$ with respect to \eqref{equivalent_H1_norm}, the functions $\hat u_j:=\sqrt{1+\sigma_j}\mathcal T(u_j)$ define a Hilbert basis of $L^2(\partial\Omega)$ with respect to its standard scalar product.
\end{thm}

We call a {\it Steklov expansion} on $\Omega$ an expression of the form
\begin{equation}\label{steklov_exp}
u=\sum_{j=1}^{\infty}a_ju_j,
\end{equation}
where $a_j:=\langle u,u_j\rangle_{H^1_{\partial}(\Omega)}$ for all $u\in H^1(\Omega)$.

From Theorem \ref{steklov_thm_0} we deduce that a Steklov expression of the form \eqref{steklov_exp} represents a function in $\mathcal H^1(\Omega)$ if and only if $\sum_{j=1}^{\infty}a_j^2<\infty$.

We recall from Theorem \ref{steklov_thm_0} that $\left\{\hat u_j\right\}_{j=1}^{\infty}$ with $\hat u_j:=\sqrt{1+\sigma_j}\mathcal T(u_j)$ is a orthonormal basis of $L^2(\partial\Omega)$. By the continuity of the trace operator $\mathcal T$ we have that
$$
\mathcal T(u)=\sum_{j=1}^{\infty}\frac{\langle u,u_j\rangle_{H^1_{\partial}(\Omega)}}{\sqrt{1+\sigma_j}}\hat u_j\,,\ \ \ \forall u\in H^1(\Omega).
$$
Hence, if $g=\mathcal T(u)$ for some $u\in H^1(\Omega)$, then $g$ has a Steklov expansion on $\partial\Omega$ analogous to \eqref{steklov_exp}, namely
\begin{equation}\label{steklov_exp0}
g=\sum_{j=1}^{\infty}g_j\hat u_j,
\end{equation}
where $g_j=\langle g,\hat u_j\rangle_{L^2(\partial\Omega)}$. Recall that a generic function $g$ defined on $\partial\Omega$ belongs to $L^2(\partial\Omega)$ if and only if $g$ can be written as in \eqref{steklov_exp0} for some $g_j\in\mathbb R$ satisfying $\sum_{j=1}^{\infty}g_j^2<\infty$. This motivates the following definition in \cite{auchmuty_steklov}.
\begin{defn}\label{HsA}
For all $s\geq 0$ we define $H^s_{A}(\partial\Omega)$ as the subspace of all functions $g\in L^2(\partial\Omega)$ with Steklov expansions as in \eqref{steklov_exp0} satisfying
$$
\sum_{j=1}^{\infty}(1+\sigma_j)^{2s}g_j^2<\infty.
$$
\end{defn}                                                           
According to Definition \ref{HsA} we define an inner product and the associated norm on $H^s_A(\partial\Omega)$:
\begin{equation*}
\langle f, g\rangle_{H^s_A(\partial\Omega)}:=\sum_{j=1}^{\infty}(1+\sigma_j)^{2s}g_jf_j\, ,
\end{equation*}
\begin{equation*}
\|g\|^2_{H^s_A(\partial\Omega)}:=\sum_{j=1}^{\infty}(1+\sigma_j)^{2s}g_j^2.
\end{equation*}
An extension operator $E:H^{1/2}_A(\partial\Omega)\rightarrow \mathcal H^1(\Omega)$ is defined in a natural way by setting
\begin{equation}\label{extension_A}
Eg=\sum_{j=1}^{\infty}\sqrt{1+\sigma_j}g_ju_j.
\end{equation}
We have the following theorem.
\begin{thm}
If $s=\frac{1}{2}$ then $  H^{s }_A(\partial\Omega)=\mathcal T(H^1(\Omega))$ 
\end{thm}
To prove the inclusion $H^{1/2}_A(\partial\Omega)\subset\mathcal T(H^1(\Omega))$ it is sufficient to show that for any $u\in H^{1/2}_A(\partial\Omega)$ there exists $u\in H^1(\Omega)$ such that $\mathcal T(u)=g$. This is done by using the extension operator \eqref{extension_A} and by setting $u=Eg\in\mathcal H^1(\Omega)$. It is standard to prove that $\mathcal T(u)=g$. Proving the reverse inclusion consists in proving that $\mathcal T(u)\in H^{1/2}_A(\partial\Omega)$ whenever $u\in H^1(\Omega)$. This is done by noting that $u=w+v$ with $w\in H^1_0(\Omega)$ and $v\in\mathcal H^1(\Omega)$. Thus $\mathcal T(u)=\mathcal T(v)$. It is then sufficient to write the Steklov expansions for $v$ and $\mathcal T(v)$ and to check the summability conditions for the Steklov coefficients. We refer to \cite{auchmuty_steklov} for a detailed proof.

In particular, from Theorem \ref{trace_thm_gagliardo} it follows that the space $H^{1/2}_A(\partial\Omega)$ coincides with the classical Sobolev space of fractional order $H^{1/2}(\partial\Omega)$. We do not know whether the spaces $H^s_A(\partial\Omega)$ provide trace spaces for higher order Sobolev spaces, in other words, if they coincide with $H^s(\partial\Omega)$ for $s\ne\frac{1}{2}$.

 In Definition \ref{HsA} the asymptotic behaviour of the eigenvalues $\sigma_j$ as $j\to \infty$ plays a crucial role.  In view of this, we recall that in the case of smooth domains, Steklov eigenvalues satisfy the following Weyl's asymptotic law
\begin{equation*}
\sigma_j\sim\frac{2\pi}{\omega_{N-1}^{\frac{1}{N-1}}}\left(\frac{j}{|\partial\Omega|}\right)^{\frac{1}{N-1}}\,,\ \ \ {\rm as\ }j\rightarrow+\infty,
\end{equation*}
where $\omega_{N-1}$ denotes the volume of the unit ball in $\mathbb R^{N-1}$. Hence we can identify the space $H^{1/2}_A(\partial\Omega)$ with the space of sequences
\begin{equation}\label{weyl2}
\left\{(s_j)_{j=1}^{\infty}\in\mathbb R^{\infty}: (j^{\frac{1}{2(N-1)}}s_j)_{j=1}^{\infty}\in l^2\right\}.
\end{equation}
Note the natural appearance of the exponent $\frac{1}{2}$ in \eqref{weyl2}.  It remarkable that, `mutatis mutandis', the summability condition in \eqref{weyl2} is already present in  \cite[Formula~(3)]{hadamard}  for the case of the unit disk of the plane.

\subsubsection{Case $k\geq 2$}\label{sub:Hm} Let $\Omega$ be a bounded domain in $\mathbb R^N$ of class $C^{0,1}$ and   $k\geq 2$ be fixed. We consider the following family of variational eigenvalue problems indexed by $\ell\in\left\{0,...,k-1\right\}$:
\begin{equation}\label{multi-Steklov}
\int_{\Omega}D^{k}u:D^k\varphi dx+\sum_{\substack{j=0,\\j\ne\ell}}^{k-1}\beta^{(\ell)}_j\int_{\partial\Omega}\frac{\partial^j u}{\partial\nu^j}\frac{\partial^j \varphi}{\partial\nu^j}d\sigma=\sigma^{(\ell)}\int_{\partial\Omega}\frac{\partial^{\ell} u}{\partial\nu^{\ell}}\frac{\partial^{\ell} \varphi}{\partial\nu^{\ell}}d\sigma\, ,
\end{equation}
$\forall\varphi\in H^k(\Omega)$, in the unknowns $u\in H^k(\Omega)$, $\sigma^{(\ell)}\in\mathbb R$, where $\beta^{(\ell)}_j>0$ are fixed constants for all $j=0,...,k-1$, $j\ne\ell$. Here $D^ku:D^k\varphi:=\sum_{|\alpha|=k}D^{\alpha}uD^{\alpha}\varphi$. For simplicity, we will set  $\beta^{(\ell)}_j=1$ for all $j,\ell=0,...,k-1$, $j\ne\ell$. All the results which we present remain valid for different positive values of $\beta^{(\ell)}_j$.
On $H^k(\Omega)$ we consider the scalar product
\begin{equation}\label{scalar_multi_steklov}
\langle u,v\rangle_{H^k_{\partial}(\Omega)}:=\int_{\Omega}D^ku:D^kvdx+\sum_{j=0}^{k-1}\int_{\partial\Omega}\frac{\partial^j u}{\partial\nu^j}\frac{\partial^j \varphi}{\partial\nu^j}d\sigma\,,\ \ \ \forall u,v\in H^k(\Omega)
\end{equation}
which induces the norm
\begin{equation}\label{norm_multi_steklov}
\|u\|_{H^k_{\partial}(\Omega)}^2:=\int_{\Omega}|D^ku|^2dx+\sum_{j=0}^{k-1}\int_{\partial\Omega}\left(\frac{\partial^j u}{\partial\nu^j}\right)^2d\sigma\,,\ \ \ \forall u\in H^k(\Omega).
\end{equation}
It is easy to see that the norm \eqref{norm_multi_steklov} is equivalent to the standard norm of $H^k(\Omega)$.

For all $\ell=0,...,k-1$ we denote by $H^k_{\ell,0}(\Omega)$ the closed subspace of $H^k(\Omega)$ defined by
$$
H^k_{0,\ell}(\Omega):=\left\{u\in H^k(\Omega):\gamma_{\ell}(u)=0\right\},
$$
and by $\mathcal H^k_{\ell}(\Omega)$ the orthogonal complement of $H^k_{0,{\ell}}(\Omega)$ in $H^k(\Omega)$ with respect to \eqref{scalar_multi_steklov}, namely
$$
\mathcal H^k_{\ell}(\Omega):=\left\{u\in H^k(\Omega):\langle u,v\rangle_{H^k_{\partial}(\Omega)}=0\,,\forall v\in H^k_{0,{\ell}}(\Omega)\right\}.
$$
We are ready to state the following theorem.

\begin{thm}\label{steklov_thm_k}
Let $\Omega$ be a bounded domain in $\mathbb R^N$ of class $C^{0,1}$. The eigenvalues of problem \eqref{multi-Steklov} have finite multiplicity and are given by a non-decreasing sequence of non-negative real numbers $\sigma^{(\ell)}_j$ defined by
\begin{equation*}
\sigma^{(\ell)}_j=\min_{\substack{U\subset H^k(\Omega)\setminus H^k_{0,\ell}(\Omega)\\{\rm dim}\,U=j}}\max_{\substack{u\in U\\u\ne 0}}\frac{\int_{\Omega}|D^ku|^2dx+\sum_{\substack{j=0\\j\ne\ell}}^{k-1}\int_{\partial\Omega}\left(\frac{\partial^j u}{\partial\nu^j}\right)^2d\sigma}{\int_{\partial\Omega}\left(\frac{\partial^{\ell}u}{\partial\nu^{\ell}}\right)^2d\sigma},
\end{equation*}
where each eigenvalue is repeated according to its multiplicity. Moreover, there exists a Hilbert basis $\{u_j^{(\ell)}\}_{j=1}^{\infty}$ of $\mathcal H^k_{\ell}(\Omega)$ of eigenfunctions $u_j^{(\ell)}$. Finally, by normalizing the eigenfunctions $u_j^{(\ell)}$ with respect to \eqref{norm_multi_steklov}, the functions $\hat u_j^{(\ell)}:=\sqrt{1+\sigma_j^{(\ell)}}\gamma_{\ell}(u_j^{(\ell)})$ define a Hilbert basis of $L^2(\partial\Omega)$ with respect to its standard scalar product.
\end{thm}
We refer to \cite{lamberti_provenzano_traces} for the proof of Theorem \ref{steklov_thm_k} in the case $k=2$. The proof for $k\geq 3$ is similar.

We note that in the case $k=2$, $\ell=0$ the first eigenvalue $\sigma_1^{(0)}=0$ has multiplicity one and the corresponding eigenfunctions are the constant functions on $\Omega$, while for $\ell=1$ the first eigenvalue $\sigma_1^{(1)}$ is positive. For $k\geq 3$ it is not straightforward to study the kernel of the operator. In fact all eigenfunctions corresponding to an eigenvalue $\sigma^{(\ell)}=0$ are of the form $u=\sum_{|\alpha|\leq k-1}a_{\alpha}x^{\alpha}$ for some $a_{\alpha}\in\mathbb R$ and satisfy $\frac{\partial^j u}{\partial\nu^j}=0$ on $\partial\Omega$ for all $j=1,...,k-1$, $j\ne\ell$. The fact that these conditions are satisfied by certain functions may depend also on $\Omega$ when $k\geq$ 3. A simple example is $\Omega=\left\{x\in\mathbb R^N: |x|<2\right\}$ and $u(x)=2-|x|^2$ which is an eigenfunction corresponding to $\sigma^{(1)}_1=0$ when $k=3$.

For all $\ell=0,...,k-1$ we define the spaces
\begin{equation}\label{multi-trace-spaces}
H^{k-\ell-\frac{1}{2}}_{\mathcal A}(\partial\Omega)=\left\{g\in L^2(\partial\Omega):g=\sum_{j=1}^{\infty}g_j\hat u_j^{(\ell)}{\rm\ such\ that\ }\sum_{j=1}^{\infty}(1+\sigma_j^{(\ell)})g_j^2<\infty\right\},
\end{equation}
which should not be confused with the spaces $H^s_A(\partial\Omega)$ in Definition \ref{HsA}. These spaces are endowed with a natural scalar product and an induced norm defined by
\begin{equation*}
\langle f, g\rangle_{H^{k-\ell-\frac{1}{2}}_{\mathcal A}(\partial\Omega)}:=\sum_{j=1}^{\infty}(1+\sigma_j^{(\ell)})g_jf_j
\end{equation*}
\begin{equation*}
\|g\|^2_{H^{k-\ell-\frac{1}{2}}_{\mathcal A}(\partial\Omega)}:=\sum_{j=1}^{\infty}(1+\sigma_j^{(\ell)})g_j^2,
\end{equation*}
and allow to describe the trace spaces of $H^k(\Omega)$. Namely, we have the following.
\begin{thm}\label{thm_single_traces}
Let $\Omega$ be a bounded domain in $\mathbb R^N$ of class $C^{0,1}$. Then
\begin{equation*}
\gamma_{\ell}(H^k(\Omega))=\gamma_{\ell}(\mathcal H^k_{\ell}(\Omega))=H^{k-\ell-\frac{1}{2}}_{\mathcal A}(\partial\Omega),
\end{equation*}
for all $\ell=0,...,k-1$. 

If $\Omega$ is of class $C^{k,1}$, then
\begin{equation*}
\Gamma(H^k(\Omega))=\prod_{\ell=0}^{k-1}H^{k-\ell-\frac{1}{2}}_{\mathcal A}(\partial\Omega),
\end{equation*}
and in particular
\begin{equation*}
H^{k-\ell-\frac{1}{2}}_{\mathcal A}(\partial\Omega)=H^{k-\ell-\frac{1}{2}}(\partial\Omega),
\end{equation*}
for all $\ell=0,...,k-1$.
\end{thm}
We note that the definition of the spaces $H^{k-\ell-1/2}_{\mathcal A}(\partial\Omega)$ require that $\Omega$ is of class $C^{0,1}$, which is a minimal assumption for the validity of  Theorem \ref{steklov_thm_k}. On the other hand, for the classical definition of the trace spaces $H^{k-\ell-1/2}(\partial\Omega)$ via Fourier analysis or Besov spaces we need $\Omega$ to be at least of class $C^{k-l-1,1}$.        We refer to \cite{lamberti_provenzano_traces} for the proof of Theorem \ref{thm_single_traces} in the case $k=2$. The proof for $k\geq 3$ can be carried out by following  the same lines.

Theorem \ref{thm_single_traces} implies that for a domain of class $C^{0,1}$
\begin{equation}\label{inclusion}
\Gamma(H^k(\Omega))\subset \prod_{\ell=0}^{k-1}H^{k-\ell-\frac{1}{2}}_{\mathcal A}(\partial\Omega).
\end{equation}
This provides a necessary condition for an element $g\in (L^2(\partial\Omega))^k$ to be the total trace $\Gamma(u)$ of some $u\in H^k(\Omega)$. This condition is not in general sufficient. In fact, as pointed out in \cite{lamberti_provenzano_traces} in the case $k=2$, the inclusion \eqref{inclusion} is in general strict, if the domain is not of class $C^{k,1}$. It has been shown in \cite{lamberti_provenzano_traces}, as one expects, that further compatibility conditions are required. These compatibility conditions can be written in a compact implicit from, as stated in the next theorem.


\begin{thm}\label{comp_thm}
Let $\Omega$ be a bounded domain in $\mathbb R^N$ of class $C^{0,1}$. Let $(g^{(0)},...,$ $ g^{(k-1)})\in\prod_{\ell=0}^{k-1}H^{k-\ell-1/2}_{\mathcal A}(\partial\Omega)$ be given by
$$
g^{(\ell)}=\sum_{j=1}^{\infty}g^{(\ell)}_j\hat u_j^{(\ell)},
$$
with $\sum_{j=1}^{\infty}(1+\sigma_j^{(\ell)})(g^{(\ell)}_j)^2<\infty$, for all $\ell=0,...,k-1$. Then $(g^{(0)},...,g^{(k-1)})$ belongs to $\Gamma(H^k(\Omega))$ if and only if for some $\ell\in\left\{0,...,k-1\right\}$
\begin{equation}\label{compatibility}
\left(\sum_{j=1}^{\infty}\sqrt{1+\sigma_j^{(\ell)}}g_j^{(\ell)}\gamma_m(u_j^{(\ell)})-g^{(m)}\right)_{m=1}^{k-1}\in \Gamma(H^{k}_{0,\ell}(\Omega)).
\end{equation}
\end{thm}
The proof is carried out by noting that, if $(g^{(0)},...,g^{(k-1)})$ belongs to the space $\Gamma(H^k(\Omega))$, then $g^{(\ell)}=\gamma_{\ell}(u_{\ell}+u_{0,\ell})$, where $u_{0,\ell}\in H^k_{0,\ell}(\Omega)$ and $u_{\ell}=\sum_{j=1}^{\infty}(1$ $ +\sigma_j^{(\ell)})g_j^{(\ell)}u_j^{(\ell)}$. We deduce that $\gamma_m(u_{\ell})-g^{(m)}=-\gamma_m(u_{0,\ell})$ for all $m=1,...,k-1$, $m\ne\ell$, and therefore the validity of \eqref{compatibility}.

We remark that it is sufficient that \eqref{compatibility} holds for just one $\ell\in\left\{0,...,k-1\right\}$. If this is true, then \eqref{compatibility} holds for all $\ell\in\left\{0,...,k-1\right\}$. Note that the problem is reduced by one dimension by condition \eqref{compatibility} because the entry corresponding to the index $m=\ell$ in the left-hand side of \eqref{compatibility} is zero.

Condition \eqref{compatibility} is quite implicit, however it is possible to re-formulate it in a more explicit, recursive way. In fact, we note that \eqref{compatibility} allows to reduce the study of $\Gamma(H^k(\Omega))$ to the study of $\Gamma(H^k_{0,\ell}(\Omega))$ for some $\ell\in\left\{0,...,k-1\right\}$. Then, we may replace through all Subsection \ref{sub:Hm} the space $H^k(\Omega)$ by $H^k_{0,\ell}(\Omega)$ and perform the same analysis. In particular, we can introduce families of polyharmonic Steklov-type problems as in \eqref{multi-Steklov}, where the variational problem is taken in $H^k_{\ell,0}(\Omega)$, and replace $\ell$ by some $\ell'\ne\ell$ in the right-hand side of the equality in \eqref{multi-Steklov}. Associated with this family of problems, we find suitable spaces defined on the boundary of $\Omega$ which allow to describe the trace spaces $\gamma_{\ell'}(H^k_{0,\ell}(\Omega))$ for all $\ell'\ne\ell$ by means of Fourier series. As in Theorem \ref{comp_thm}, a description of $\Gamma(H^k_{0,\ell}(\Omega))$ is deduced from the knowledge of $\gamma_{\ell'}(H^k_{0,\ell}(\Omega))$ and $\Gamma(H^k_{0,\ell}(\Omega)\cap H^k_{0,\ell'}(\Omega))$. 
Thus, the problem is reduced again by one dimension, namely, it is reduced to the study of $\Gamma(H^k_{0,\ell}(\Omega)\cap H^k_{0,\ell'}(\Omega))$. This process stops after $k-1$ steps.

To clarify the ideas, we will briefly describe the case $k=2$, for which necessary and sufficient conditions are deduced.


For $\ell=0,1$, we denote by $\mathcal B^2_{0,\ell}(\Omega)$ the orthogonal complement of $H^2_0(\Omega)=H^2_{0,0}(\Omega)\cap H^2_{0,1}(\Omega)$ in $H^2_{0,\ell}(\Omega)$ with respect to the quadratic form \eqref{scalar_multi_steklov} (with $k=2$), namely
$$
\mathcal B^2_{0,\ell}(\Omega):=\left\{u\in H^2_{0,\ell}(\Omega):\langle u,\varphi\rangle_{H^2_{\partial}(\Omega)}=0\,,\forall \varphi\in H^2_0(\Omega)\right\}.
$$
Note that $\mathcal B^2_{0,\ell}(\Omega)$ is the subspace of the  biharmonic functions in $H^2(\Omega)$ with  $\gamma_{\ell}(u)=0$. Recall that the biharmonic functions in $H^2(\Omega)$ are defined as those functions $u\in H^2(\Omega)$ such that $\int_{\Omega}D^2u:D^2\varphi dx=0$ for all $\varphi\in H^2_0(\Omega)$. Therefore, we have
$$
H^2_{0,\ell}(\Omega)=H^2_0(\Omega)\oplus\mathcal B^2_{0,\ell}(\Omega),
$$
where the sum is orthogonal with respect to \eqref{scalar_multi_steklov}. Next, we consider  the following family of auxiliary variational eigenvalue problems indexed by $\ell,m\in\left\{0,1\right\}$, $\ell\ne m$:
\begin{equation}\label{auxiliary_problem}
\int_{\Omega}D^{2}w:D^2\varphi dx=\eta^{\ell,m}\int_{\partial\Omega}\frac{\partial^{m} w}{\partial\nu^{m}}\frac{\partial^{m} \varphi}{\partial\nu^{m}}d\sigma\,,\ \ \ \forall\varphi\in H^2_{0,\ell}(\Omega),\end{equation}
in the unknowns $w\in H^2_{0,\ell}(\Omega)$ and $\eta^{\ell,m}\in\mathbb R$. We have the following theorem.
\begin{thm}\label{auxiliary_thm_k}
Let $\Omega$ be a bounded domain in $\mathbb R^N$ of class $C^{0,1}$, $\ell,m\in\left\{0,1\right\}$, $\ell\ne m$. The eigenvalues of problem \eqref{auxiliary_problem} have finite multiplicity and are given by a non-decreasing sequence of non-negative real numbers $\eta^{\ell,m}_j$ defined by
\begin{equation*}
\eta^{\ell,m}_j=\min_{\substack{W\subset H^2_{0,\ell}(\Omega)\setminus H^2_0(\Omega)\\{\rm dim}\,W=j}}\max_{\substack{w\in W\\w\ne 0}}\frac{\int_{\Omega}|D^2w|^2dx}{\int_{\partial\Omega}\left(\frac{\partial^{m}w}{\partial\nu^{m}}\right)^2d\sigma},
\end{equation*}
where each eigenvalue is repeated according to its multiplicity. Moreover, there exists a Hilbert basis $\{w_j^{\ell,m}\}_{j=1}^{\infty}$ of $\mathcal B^2_{0,\ell}(\Omega)$ of eigenfunctions $w_j^{\ell,m}$. Finally, by normalizing the eigenfunctions $w_j^{\ell,m}$ with respect to \eqref{norm_multi_steklov}, the functions $\hat w_j^{\ell,m}:=\sqrt{1+\eta_j^{\ell,m}}\gamma_m(w_j^{\ell,m})$ define a Hilbert basis of $L^2(\partial\Omega)$ with respect to its standard scalar product.
\end{thm}
We refer to \cite{lamberti_provenzano_traces} for the proof of Theorem \ref{auxiliary_thm_k}.

In order to characterize those couples $(g^{(0)},g^{(1)})\in H^{\frac{3}{2}}_{\mathcal A}(\partial\Omega)\times H^{\frac{1}{2}}_{\mathcal A}(\partial\Omega)$ whi\-ch belong to $\Gamma(H^2(\Omega))$, we need to introduce the spaces $H^{\frac{3}{2}-m}_{\mathcal A,\ell}(\partial\Omega)$ defined for  $\ell,m\in\left\{0,1\right\}$, $\ell\ne m$, by
\begin{equation*}
H^{\frac{3}{2}-m}_{\ell,\mathcal A}(\partial\Omega):=\left\{g\in L^2(\partial\Omega):g=\sum_{j=1}^{\infty}g_j\hat w_j^{\ell,m}{\rm\ such\ that\ }\sum_{j=1}^{\infty}(1+\eta_j^{\ell,m})g_j^2<\infty\right\},
\end{equation*}
where $g_j=\langle g,\hat w_j^{\ell,m}\rangle_{L^2(\partial\Omega)}$. It turns out that $H^{3/2-m}_{\ell,\mathcal A}(\partial\Omega)=\gamma_{m}(\mathcal B^2_{0,\ell}(\Omega))$.

We have the following corollary of Theorem \ref{comp_thm}

\begin{cor}\label{comp_thm_2}
Let $\Omega$ be a bounded domain in $\mathbb R^N$ of class $C^{0,1}$. Let $(g^{(0)},g^{(1)})$ $ \in H^{3/2}_{\mathcal A}(\partial\Omega)\times H^{1/2}_{\mathcal A}(\partial\Omega)$ be given by
$$
g^{(\ell)}=\sum_{j=1}^{\infty}g^{(\ell)}_j\hat u_j^{(\ell)},
$$
with $\sum_{j=1}^{\infty}(1+\sigma_j^{(\ell)})(g^{(\ell)}_j)^2<\infty$, for $\ell=0,1$. Then $(g^{(0)},g^{(1)})$ belongs to the space $\Gamma(H^2(\Omega))$ if and only if one of the following two equivalent conditions holds:
\begin{equation}\label{compatibility2}
\sum_{j=1}^{\infty}\sqrt{1+\sigma_j^{(0)}}g_j^{(0)}\gamma_1(u_j^{(0)})-g^{(1)}\in H^{\frac{1}{2}}_{\mathcal A,0}(\partial\Omega).
\end{equation}
\begin{equation}\label{compatibility23}
\sum_{j=1}^{\infty}\sqrt{1+\sigma_j^{(1)}}g_j^{(1)}\gamma_0(u_j^{(1)})-g^{(0)}\in H^{\frac{3}{2}}_{\mathcal A,1}(\partial\Omega).
\end{equation}
\end{cor}

We conclude this section with some remarks.  We note that for $k=2$, in the case of smooth domains, the eigenvalues $\sigma_j^{(0)}$ and $\sigma_j^{(1)}$ satisfy the following asymptotic law
$$
\sigma_j^{(0)}\sim C_N\left(\frac{j}{|\partial\Omega|}\right)^\frac{3}{N-1}{\rm\ \ \ and\ \ \ }\sigma_j^{(1)}\sim C_N'\left(\frac{j}{|\partial\Omega|}\right)^\frac{1}{N-1}\,,{\rm\ as\ }j\rightarrow+\infty,
$$
where $C_N,C_N'$ depend only on $N$, see \cite{lamberti_provenzano_traces} for details. Hence we can identify the space $H^{3/2}_{\mathcal A}(\partial\Omega)$ with the space of sequences
\begin{equation}\label{weyl1}
\left\{(s_j)_{j=1}^{\infty}\in\mathbb R^{\infty}:(j^{\frac{3}{2(N-1)}}s_j)_{j=1}^{\infty}\in l^2\right\}
\end{equation}
and the space $H^{1/2}_{\mathcal A}(\partial\Omega)$ with
\begin{equation}\label{weyl20}
\left\{(s_j)_{j=1}^{\infty}\in\mathbb R^{\infty}:(j^{\frac{1}{2(N-1)}}s_j)_{j=1}^{\infty}\in l^2\right\}.
\end{equation}
Also in this case we observe the natural  appearance of the exponents $\frac{3}{2}$ and $\frac{1}{2}$ in \eqref{weyl1} and \eqref{weyl20} respectively.

Note that for $k=2$ and $\ell =0$, setting $\lambda:=-\beta^{(0)}_1$, problem \eqref{multi-Steklov} is the weak formulation of the following Steklov-type problem for the biharmonic operator:
\begin{equation}\label{classic_biharmonic_steklov_1}
\begin{cases}
\Delta^2 u=0, & {\rm in\ }\Omega,\\
\frac{\partial^2 u}{\partial\nu^2}-\lambda\frac{\partial u}{\partial\nu}=0, & {\rm on\ }\partial\Omega,\\
-{\rm div}_{\partial\Omega}(D^2u\cdot\nu)_{\partial\Omega}-\frac{\partial\Delta u}{\partial\nu}=\sigma^{(0)}(\lambda) u, & {\rm on\ }\partial\Omega,
\end{cases}
\end{equation}
in the unknowns $u$ (the eigenfunction) and $\sigma^{(0)}(\lambda)$ (the eigenvalue). Similarly,  for $k=2$ and $\ell =1$,    setting  $\mu:=-\beta^{(1)}_0$, problem \eqref{multi-Steklov} is the weak formulation of the following Steklov-type problems for the biharmonic operator:
\begin{equation}\label{classic_biharmonic_steklov_2}
\begin{cases}
\Delta^2 u=0, & {\rm in\ }\Omega,\\
\frac{\partial^2 u}{\partial\nu^2}=\sigma^{(1)}(\mu)\frac{\partial u}{\partial\nu}, & {\rm on\ }\partial\Omega,\\
-{\rm div}_{\partial\Omega}(D^2u\cdot\nu)_{\partial\Omega}-\frac{\partial\Delta u}{\partial\nu}-\mu u=0, & {\rm on\ }\partial\Omega,
\end{cases}
\end{equation}
 in the unknowns $u$ (the eigenfunction) and $\sigma^{(1)}(\mu)$ (the eigenvalue).  According to \eqref{multi-Steklov}, the numbers $\lambda,\mu$ are assumed to be strictly negative. Problems \eqref{classic_biharmonic_steklov_1} and \eqref{classic_biharmonic_steklov_2} admit increasing sequences of eigenvalues, which we denote by $\{\sigma^{(0)}_j(\lambda)\}_{j=1}^{\infty}$ and  $\{\sigma^{(1)}_j(\mu)\}_{j=1}^{\infty}$ respectively. We have highlighted the dependence of the eigenvalues on $\lambda$ and $\mu$. Indeed, problems \eqref{classic_biharmonic_steklov_1} and \eqref{classic_biharmonic_steklov_2} define a family of multi-parameter Steklov-type problems for the biharmonic operator, which are  genuine generalizations of the classical Steklov problem \eqref{Steklov_classic} for the Laplace operator.

As pointed out in \cite{lamberti_provenzano_traces}, problems \eqref{classic_biharmonic_steklov_1} and \eqref{classic_biharmonic_steklov_2} can be studied also for $\lambda,\mu\geq 0$. We refer to \cite{lamberti_provenzano_traces} for a detailed analysis of the dependence of the eigenvalues $\sigma^{(0)}_j(\lambda)$ and $\sigma^{(1)}_j(\mu)$ upon $\lambda,\mu$, explicit examples, Weyl's asymptotics, as well as the asymptotic behavior of the eigenvalues for $\lambda,\mu\rightarrow-\infty$.

For $k=2$, the auxiliary problems defined by \eqref{auxiliary_problem} are the weak formulations of the following Steklov-type problems for the biharmonic operator. For $\ell=1$, $m=0$ we have
\begin{equation}\label{NBS}
\begin{cases}
\Delta^2 w=0, & {\rm in\ }\Omega,\\
\frac{\partial w}{\partial\nu}=0, & {\rm on\ }\partial\Omega,\\
-{\rm div}_{\partial\Omega}(D^2w\cdot\nu)_{\partial\Omega}-\frac{\partial\Delta w}{\partial\nu}=\eta^{1,0} w, & {\rm on\ }\partial\Omega,
\end{cases}
\end{equation}
in the unknowns $w$ (the eigenfunction) and $\eta^{1,0}$ (the eigenvalue). For $\ell=0$, $m=1$ we have 
\begin{equation}\label{DBS}
\begin{cases}
\Delta^2 w=0, & {\rm in\ }\Omega,\\
\frac{\partial^2 w}{\partial\nu^2}=\eta^{0,1}\frac{\partial w}{\partial\nu}, & {\rm on\ }\partial\Omega,\\
w=0, & {\rm on\ }\partial\Omega,
\end{cases}
\end{equation}
in the unknowns $w$ (the eigenfunction) and $\eta^{0,1}$ (the eigenvalue).

Problem \eqref{classic_biharmonic_steklov_1} for $\lambda =0$ has been introduced in \cite{buopro} and further investigated in \cite{buoprocha}. 
Problem \eqref{DBS} has been considered by many authors in the literature, while problem \eqref{NBS} has been much less investigated. See \cite{lamberti_provenzano_traces} for references, see also \cite{buoso}. In particular, it has been proved in \cite{lamberti_provenzano_traces} that problem \eqref{NBS} is the limit problem of \eqref{classic_biharmonic_steklov_1} as $\lambda\rightarrow-\infty$, while problem \eqref{DBS} is the limit problem of \eqref{classic_biharmonic_steklov_2} as $\mu\rightarrow-\infty$.

Already for $k=3$,  writing explicitly the classical formulation of problem \eqref{multi-Steklov} analogous to \eqref{classic_biharmonic_steklov_1}-\eqref{classic_biharmonic_steklov_2}  is not easy, see e.g., \cite{arfela, lamberti_ferraresso} for a discussion of various classical boundary value problems for polyharmonic operators. However, even if their explicit form is quite involved, there exist uniquely defined boundary differential operators $\mathcal N_j$, $j=0,...,k-1$ of order $j+k$ such that any smooth solution to \eqref{multi-Steklov} solves the following boundary value problem
\begin{equation*}
\begin{cases}
\Delta^ku=0, & {\rm in\ }\Omega,\\
\mathcal N_{k-1-j}u+\beta_{j}^{(\ell)}\frac{\partial^j u}{\partial\nu^{j}}=0\ \forall j=0,...,k-1,j\ne\ell, & {\rm on\ }\partial\Omega,\\
\mathcal N_{k-1-\ell}u=\sigma^{(\ell)}\frac{\partial^{\ell}u}{\partial\nu^{\ell}}, & {\rm on\ }\partial\Omega,
\end{cases}
\end{equation*}
in the unknowns $u$ (the eigenfunction) and $\sigma^{(\ell)}=\sigma^{(\ell)}(\beta^{(\ell)}_0,...,,\beta^{(\ell)}_{\ell-1},\beta^{(\ell)}_{\ell+1},...,$ $ \beta^{(\ell)}_{k-1})$ (the eigenvalue). 
When $k=2$, $\mathcal N_0u=\frac{\partial^2u}{\partial\nu^2}$ and $\mathcal N_1(u)=-{\rm div}_{\partial\Omega}(D^2u\cdot\nu)_{\partial\Omega}-\frac{\partial\Delta u}{\partial\nu}$. The operators $\mathcal N_j$ correspond to the Neumann boundary conditions for the polyharmonic operator $(-\Delta)^k$.

\section*{Acknowledgements}
The authors are very thankful to Professors  Giles Auchmuty and  Victor I. Burenkov  for useful discussions and references.

\end{document}